\documentclass[english]{amsart}

\usepackage{babel}
\usepackage{amstext}
\usepackage{amsmath}
\usepackage{amsfonts}
\usepackage{latexsym}
\usepackage{ifthen}
\usepackage[all,2cell,dvips]{xy} \UseAllTwocells \SilentMatrices
\usepackage{xy}
\xyoption{all}
\newcommand{\id}{{\rm id}}

\newcommand{\codim}{{\rm codim}}

\newcommand{\im}{{\rm im}}

\newtheorem{lemma1}[equation]{}

\newenvironment{lemma}{\begin{lemma1}{\bf Lemma.}}{\end{lemma1}}
\newenvironment{example}{\begin{lemma1}{\bf Example.}\rm}{\end{lemma1}}

\newenvironment{theorem}{\begin{lemma1}{\bf Theorem.}}{\end{lemma1}}
\newenvironment{proposition}{\begin{lemma1}{\bf Proposition.}}{\end{lemma1}}
\newenvironment{corollary}{\begin{lemma1}{\bf Corollary.}}{\end{lemma1}}

\newenvironment{conjecture}{\begin {lemma1}{\bf Conjecture.}}{\end{lemma1}}

\newcommand{\Q}{\ensuremath{\mathbb{Q}}}
\newcommand{\Z}{\ensuremath{\mathbb{Z}}}
\newcommand{\C}{\ensuremath{\mathbb{C}}}

\newcommand{\PP}{\ensuremath{\mathbb{P}}}
\newcommand{\D}{\ensuremath{\mathbb{D}}}

\newcommand{\Homsheaf} { \ensuremath{ \mathcal{H} \! om}}
\newcommand{\Ssheaf} { \ensuremath{ \mathcal{S} } }
\newcommand{\Isheaf} { \ensuremath{ \mathcal{I} } }

\newcommand{\Osheaf} { \ensuremath{ \mathcal{O} } }

\newcommand{\chow}[1]{\ensuremath{\mathcal{C}(#1)}}

\newcommand{\merom}[3]{\ensuremath{#1:#2 \dashrightarrow #3}}

\newcommand{\holom}[3]{\ensuremath{#1:#2  \rightarrow #3}}
\newcommand{\fibre}[2]{\ensuremath{#1^{-1} (#2)}}

\makeatletter
\ifnum\@ptsize=0  \fi \ifnum\@ptsize=2  \fi 

\newcommand\sO{{\mathcal O}}

\DeclareMathOperator*{\Aut}{Aut}

\DeclareMathOperator*{\rk}{rk}
\DeclareMathOperator*{\red}{red}
\DeclareMathOperator*{\nons}{nons}

\newcommand{\Bl}[2]{\ensuremath{\mbox{Bl}_{#1}{#2}}}

\title {Uniruled varieties with split tangent bundle} 
\date{\today}

\begin{document}

%\begin{abstract}
%do it later
%\end{abstract}

\author{Andreas H\"oring}
\address{Andreas H\"oring, c/o Lehrstuhl Mathematik I, Mathematisches Institut,
Universit\" at Bayreuth,
D-95440 Bayreuth, Germany}
\email{andreas.hoering@uni-bayreuth.de}

\maketitle

\tableofcontents

\section{Introduction}

The goal of this paper is to give a partial answer to the following conjecture, initially asked by A. Beauville.

\begin{conjecture}
\label{conjecturebeauville}
Let $X$ be a compact K\"ahler manifold such that $T_X = V_1 \oplus V_2$, where $V_1$ and $V_2$ are vector bundles.
Let \holom{\mu}{\tilde{X}}{X} be the universal covering of $X$. Then $\tilde{X} \simeq Y_1 \times Y_2$, where
$\dim Y_j = \rk V_j$. If moreover $V_j$ is integrable,  $\mu^* V_j \simeq p_{Y_j}^* T_{Y_j}$ (up to an appropriate
automorphism of $\tilde{X}$).
\end{conjecture}

The integrability of the vector bundles $V_j$ is not true in general (see example \ref{examplebeauville}),
but the decomposition of the universal covering $\tilde{X}$ still holds in these examples. The conjecture
has been studied before by Beauville \cite{Bea00}, Druel \cite{Dru00}, Campana-Peternell \cite{CP02}
and  recently by Brunella-Pereira-Touzet \cite{BPT04}. Their paper contains most of the preceeding
results, its main result being the

\begin{theorem} \cite[Thm.1]{BPT04}
\label{theorembrunella}
Let $X$ be a compact K\"ahler manifold. Suppose that its tangent bundle splits as $T_X = V_1 \oplus V_2$,
where $V_2 \subset T_X$ is a subbundle of rank $\dim X - 1$. There are two cases:
\begin{enumerate}
 \item if $V_2$ is not integrable,  $V_1$ is tangent to the fibres of a $\PP^1$-bundle.
 \item if $V_2$ is integrable,  \holom{\mu}{\tilde{X}}{X}, the universal covering of $X$, splits as $Y_1 \times Y_2$,
 where $Y_2$ is a curve. Moreover $\mu^* V_1 \simeq p_{Y_1}^* T_{Y_1}$ and $\mu^* V_2 \simeq p_{Y_2}^* T_{Y_2}$. 
 \end{enumerate}
\end{theorem}

The theorem establishes a surprising link between the existence of rational curves along the foliation $V_1$ and
the integrability of the complement $V_2$. We obtain a similar statement in the projective case.

\begin{theorem}
\label{theoremnonuniruledintegrable}
Let $X$ be a projective manifold with split tangent bundle $T_X=V_1 \oplus V_2$.
If $X$ is not uniruled, $V_1$ and $V_2$ are integrable.
\end{theorem} 

In this paper we will  concentrate on projective 
uniruled varieties. By the theorem above this is
the class of varieties where the integrability of the direct factors $V_j$ might fail, but we will obtain
integrability results in special cases.
The global strategy is to construct a fibre space structure on $X$ which is
related to the decomposition $T_X = V_1 \oplus V_2$. 
The classical Ehresmann theorem \ref{theoremehresmann}
 allows us to deduce 
properties of the universal covering of $X$ from this 
fibre space structure.

The structure of the paper is as follows: in the next section we recall some results about
foliations and rational curves. We also give the short proof of theorem \ref{theoremnonuniruledintegrable}.
Section \ref{sectionrationallyconnected} treats rationally connected manifolds, which are known to
be simply connected, where we obtain a rather
satisfactory answer to the conjecture.

\begin{theorem}
\label{theorembogomolovbetter}
Let $X$ be a rationally connected manifold such that $T_X = V_1 \oplus V_2$. 
If $V_1$ or $V_2$ is integrable,  
$X$ is isomorphic to a product $Y_1 \times Y_2$ such that
$V_j = p_{Y_j}^* T_{Y_j}$ for $j=1,2$.
\end{theorem}

We then study the structure of elementary Mori contractions of projective uniruled fourfolds
with split tangent bundle. Combining the results of section \ref{sectionuniruled}, we obtain

\begin{theorem}
\label{theoremdim4uniruled}
Let $X$ be a projective uniruled fourfold such that $T_X = V_1 \oplus V_2$ where $\rk V_1=\rk V_2=2$.
If $V_1$ and $V_2$ are integrable, conjecture \ref{conjecturebeauville} holds.
\end{theorem}

{\bf Acknowledgements.}
This paper is part of the author's binational Ph.D. thesis supervised by L. Bonavero (Institut Fourier Grenoble) and
T. Peternell (Universit\"at Bayreuth) which is partially supported by the DFG-Schwerpunkt \lq Globale Methoden in der Komplexen
Geometrie\rq. I want 
\nopagebreak 
to thank St\'ephane Druel for his careful proofreading of a rather awkward first version of this
paper.

\section{Notation and basic results}
\label{sectionbasics}

We work over the complex field $\C$. For standard definitions in complex algebraic geometry we refer to  \cite{Ha77} or \cite{Kau83}, 
for positivity notions of vector bundles we follow the definitions from \cite{Laz04II}. Manifolds and varieties are always supposed
to be connected.

A fibration is a surjective projective morphism \holom{\phi}{X}{Y} with connected fibres
from a quasi-projective manifold 
to a normal quasi-projective variety $Y$ 
such that $\dim X > \dim Y$. The $\phi$-smooth locus is the largest Zariski open subset $Y^* \subset Y$
such that for every $y \in Y^*$, the fibre $\fibre{\phi}{y}$ is a smooth variety of dimension $\dim X - \dim Y$.
The $\phi$-singular locus is its complement.
A fibre is always a fibre in the scheme-theoretic sense, a set-theoretic fibre is the reduction of the fibre.
We say that a fibration is almost smooth if it is equidimensional and every set-theoretic fibre is a
smooth variety.

Let \holom{\phi}{X}{Y} be a morphism from a projective manifold to a normal projective variety $Y$. A line bundle $L$ on $X$
is $\phi$-trivial if for every $\phi$-fibre $F$, we have $L|_F \simeq \sO_F$.

\subsection{Rational curves and split tangent bundle}
\label{subsectionrationalcurves}

We will use the standard terminology from Mori theory. A Mori contraction of a projective manifold $X$ is
a morphism with connected fibres \holom{\phi}{X}{Y} to a normal variety $Y$ such that the anticanonical bundle $-K_X$ is $\phi$-ample.
We say that the contraction is elementary if there exists a rational curve $C \subset X$ such that a curve
$C' \subset X$ is contracted by $\phi$ if and only if we have an equality in $N_1(X)$
\[
[C'] = \lambda \ [C] \qquad \lambda \in \Q^+.
\]
By $N_1(X)$ we denote the $\Q$-vector space of 1-cycles on $X$ modulo numerical equivalence (cf. \cite[1.3]{Deb01}).
The contraction is said to be of fibre type if $\dim Y < \dim X$; otherwise it is birational. 

\begin{proposition}
\cite[Prop.1.8.]{CP02}
\label{propositionCP02Prop18}
Let $X$ be a projective manifold such that $T_X = V_1 \oplus V_2$. 
Let \holom{\phi}{X}{Y} be an elementary 
extremal contraction and let $F$ be a fibre of $\phi$.
\begin{enumerate}
\item If $L_j = \det V_j$ is not $\phi$-trivial, we have $\dim F \leq \rk V_j$.

\item $\dim F \leq  \max (\rk V_1, \rk V_2)$.

\item Suppose $F$ contains a rational curve \holom{f}{\PP^1}{X} such that
\begin{equation}
\label{equationsplittingtype}
f^* T_X \simeq \Osheaf_{\PP^1} (2) \oplus \bigoplus_{i=2}^{\dim X} \Osheaf_{\PP^1} (a_i)  
\qquad \mbox{with} \quad a_i \leq 1 \qquad \forall \ i>1.
\end{equation}
Then, up to renumbering, $L_1$ is $\phi$-ample, $L_2$ is $\phi$-trivial, and $\Osheaf_{\PP^1} (2) \subset f^* V_1$.
\end{enumerate}
\end{proposition}

{\bf Remark.} If \holom{\phi}{X}{Y} is an elementary contraction of fibre type, 
a general fibre $F$ always contains a rational curve of splitting type (\ref{equationsplittingtype}).
Indeed $F$ is a Fano manifold of positive dimension, so by
\cite[Ex.4.8.3]{Deb01} there exists a rational curve \holom{f}{\PP^1}{F\subset X} such that
$f^* T_F$ is nef and has splitting type (\ref{equationsplittingtype}).
Conclude with lemma \ref{lemmasplit}.

\begin{lemma}
\label{lemmasplit}
Let \holom{\phi}{X}{Y} be a morphism from a projective manifold $X$
to a normal variety $Y$.  
Let $F$ be a smooth irreducible component of a reduced $\phi$-fibre $Z$. 
Let \holom{f'}{\PP^1}{F} be a rational curve 
such that $f'^* T_F$ is nef. There exists a deformation \holom{f}{\PP^1}{F} of $f'$ in $F$, such that  
$f^* T_F$ is nef and 
\begin{equation}
\label{equationfibresplit}
f^* T_X \simeq f^* T_F \oplus f^* N_{F/X}.
\end{equation}
\end{lemma}

{\bf Proof. }
Since $Z$ is reduced, the canonical morphism
\[
\Isheaf_Z/\Isheaf_Z^2 \otimes \Osheaf_F \rightarrow \Isheaf_F/\Isheaf_F^2
\]
is generically surjective. Since $\Isheaf_Z/\Isheaf_Z^2$ is globally generated, 
$\Isheaf_F/\Isheaf_F^2=N_{F/X}^*$ is generically generated.
Since $F$ is smooth, the deformations of $f'$ cover $F$ (\cite[Prop.4.8]{Deb01}). So for 
a general deformation \holom{f}{\PP^1}{F}, 
the restriction of the bundle $N_{F/X}^*$ to $f(\PP^ 1)$ is generically generated
and $f^* T_F$ is nef. It follows that $f^* (N_{F/X}^* \otimes T_F)$ is nef,
so $H^1(\PP^1, f^* (N_{F/X}^* \otimes T_F))=0$. Hence
the exact sequence
\[
0 \rightarrow f^* T_F  \rightarrow f^* T_X \rightarrow f^* N_{F/X} \rightarrow  0
\]
splits. $\square$

\subsection{Foliations}
\label{subsectionfoliations}

In this section we introduce the terminology that we will need and state some  results
about holomorphic foliations. Since the analytic category provides the adapted framework for this theory 
we state the results for compact K\"ahler manifolds.

Let $X$ be a compact K\"ahler manifold. 
A subbundle $V \subset T_X$ is integrable if it is closed under the Lie bracket.
We recall that the Lie bracket
\[
[.,.] \ : \ V  \times V \ \rightarrow \ T_X
\]
is a bilinear antisymmetric mapping that is not $\Osheaf_X$-linear but 
induces an $\Osheaf_X$-linear map $\wedge^2 V \rightarrow T_X / V$ that is zero if and only if $V$ is integrable.
In particular 
\[
H^0(X, \Homsheaf( \wedge^2 V, T_X / V)) = 0
\] 
implies that $V$ is integrable. In general we will show this vanishing property using a covering family $S$ 
(cf. \cite[Defn.1.8.]{Ca04} for the terminology of covering families) of
subvarieties $(Z_s)_{s \in S}$ of $X$ such that a general member of the family satisfies
\[
H^0(Z_s, \Homsheaf( \wedge^2 V, T_X / V)|_{Z_s}) = 0.
\] 
By the Frobenius theorem an integrable subbundle $V$ of $T_X$ induces a foliation on $X$, i.e. for every $x \in X$ there exists
an analytic neighbourhood $U$ and a submersion \holom{q}{U}{W} such that $T_{U/W}=V|_U$. This shows that $V$ can be realised
as the tangent bundle of locally closed subsets of $X$. The foliation induces an equivalence relation
on $X$, two points being equivalent if and only if they can be connected by chains of smooth (open) curves $C_i$
such that $T_{C_i} \subset V|_{C_i}$. An equivalence class is  called a leaf of the foliation. 
The equivalence relation induces a quotient map $X \rightarrow X/V$ onto
the so-called leaf space $X/V$ such that the fibres are the leaves of the foliation. This map can
always be defined topologically, but in general $X/V$ is not Hausdorff.
A subset $X^* \subset X$ is saturated if every leaf of the foliation is either contained in $X^*$ or disjoint from it.
We will say that the general leaf of a foliation is compact if there exists a non-empty saturated Zariski open subset
$X^* \subset X$ such that every leaf contained in $X^*$ is compact.

\begin{proposition}
\label{propositiongeneralleafcompact}
Let $X$ be a compact K\"ahler manifold, and let $V \subset T_X$ be an integrable subbundle. Assume that the general $V$-leaf
is compact. Then every $V$-leaf of the foliation is compact. There exists an almost smooth holomorphic map $X \rightarrow Y$
whose set-theoretic fibres are $V$-leaves.
\end{proposition}

{\bf Proof.}
The compactness of the leaves follows from the global stability theorem on K\"ahler manifolds (cf. \cite[Thm.1]{Pe01}
for a short proof). 
Holmann \cite[Cor.3.2]{Ho80} has shown that in this case 
the leaf space $X/V$ admits the structure of an analytic space such that the projection is almost smooth. 
$\square$

\begin{corollary}
\label{corollaryquotientsubmersion}
Let $X$ be a compact K\"ahler manifold, and let $V \subset T_X$ be an integrable subbundle such that one leaf
is compact and rationally connected. Then there exists a submersion $X \rightarrow Y$ onto a smooth variety $Y$ 
such that $T_{X/Y}=V$.
\end{corollary}

{\bf Remark.}
Kebekus, Sol\'a and Toma obtained independently similar results (\cite[Thm.28]{KST05}) for singular foliations on projective manifolds.

{\bf Proof.}
A rationally connected leaf is simply connected, so by Reeb's
local stability theorem (cf. \cite[IV, \S 3, Thm.3]{CLN85}), the general $V$-leaf
is compact. Proposition \ref{propositiongeneralleafcompact} yields an almost
smooth map \holom{g}{X}{Y} onto the leaf space $Y:=X/V$. 

Let $y \in Y$ be an arbitrary point. By \cite[Prop.3.7.]{Mol88}
there exists a neighbourhood $U$ of $y$ which is isomorphic to $\D/G$ where $\D$ is the $\dim Y$-dimensional
unit disc and $G$ is a finite group. More precisely $G$ is the image of a representation
of the fundamental group $\pi_1(\fibre{g}{y}_{\red})$.  Denote now by
$\holom{q}{\D}{U}$ the quotient map and by $X_U$ the normalisation of $\fibre{g}{U} \times_U \D$. 
Let \holom{g'}{X_U}{\D} be the induced map and \holom{q'}{X_U}{\fibre{g}{U}} the induced \'etale
covering. Then $g'$ is a submersion and the $g'$-fibres are the
$q'^* V$-leaves.

A $V$-leaf is rationally connected if and only if the
corresponding $q'^* V$-leaf is rationally connected. Furthermore
one $g'$-fibre is rationally connected if and only if all
$g'$-fibres are rationally connected (by deformation invariance of
rational connectedness, cf. \cite[Thm.3.11]{Ko96}). Now by hypothesis there
exists one rationally connected $V$-leaf, so by connectedness of
$Y$ all $V$-leaves are rationally connected. Since rationally
connected varieties are simply connected, the group $G$ is trivial.
It follows that $U\simeq\D/G=\D$ is smooth and the quotient map $\holom{q}{\D}{U}$ is an isomorphism.
Hence \holom{q'}{X_U}{\fibre{g}{U}} is also an isomorphism, so $g'$ and $g|_{\fibre{g}{U}}$ identify under
this isomorphism. In particular  $g|_{\fibre{g}{U}}$ is a submersion. Since $y$ is arbitrary, this shows the claim.
$\square$

\begin{corollary}
\label{corollaryquotientcodim0}
Let $X$ be a compact K\"ahler manifold and $V \subset T_X$ a subbundle. Suppose that
there exists a fibration \holom{\phi'}{X}{Y'} such that $T_F=V|_F$ for a general fibre $F$. Then $V$ is integrable
with compact leaves. The almost smooth map \holom{\phi}{X}{Y} from proposition 
\ref{propositiongeneralleafcompact} is a factorisation
of $\phi'$, i.e. there exists a birational morphism \holom{g}{Y}{Y'} such that $g \circ \phi = \phi'$. 
In particular if $\phi'$ is equidimensional, it is almost smooth. $\square$
\end{corollary}

\subsection{An integrability result}
\label{subsectionintegrability}

We show theorem \ref{theoremnonuniruledintegrable} and give a counterexample in the uniruled case. This shows that 
the restriction to the uniruled case is appropriate to the nature of the problem.

\begin{theorem} \cite[Thm.1.5]{CP04}
Let $X$ be a projective manifold with split tangent bundle $T_X=V_1 \oplus V_2$.
Set $L_j:= \det V_j^*$ for $j=1,2$. If $X$ is not uniruled,  $L_j$ is pseudoeffective.
\end{theorem}

{\bf Proof of theorem \ref{theoremnonuniruledintegrable}.}
By the theorem of Campana-Peternell above, the line bundle $L_1:=\det V_1^*$ is pseudoeffective.
Since $V_1^*$ is a direct factor of $\Omega_X$, the vector bundle $\det V_1 \otimes \wedge^{\rk V_1} \Omega_X$
has a trivial direct factor. 
If $\theta \in H^0(X, L_1^{-1} \otimes \wedge^{\rk V_1} \Omega_X)$ is the associated nowhere-vanishing
$\det V_1$ -valued form, and $\zeta$ a germ of any vector field, a local computation shows that $i_\zeta \theta = 0$
if and only if $\zeta$ is in $V_2$. An integrability criterion by Demailly \cite[Thm.]{De02} shows that $V_2$ is integrable.
The statement for $V_1$ follows completely analogously.
$\square$

{\bf Remark.}
The integrability theorem \ref{theoremnonuniruledintegrable} is optimal, in the sense that there is 
the following counterexample to the integrability in the uniruled case. 

\begin{example} (A. Beauville)
\label{examplebeauville}
Let $A$ be an abelian surface and  $u_1, u_2$ be linearly independent vector fields on $A$. Let $z_1, z_2$ 
be nonzero vector fields on $\PP^1$ such that $[z_1,z_2] \neq 0$.
Then $v_1:=p_A^*(u_1) + p_{\PP_1}^* (z_1)$ and $v_2:=p_A^*(u_2) + p_{\PP_1}^* (z_2)$ are 
everywhere nonzero, linearly independent 
vector fields
on $X:=A \times \PP^1$. The subbundle $V:= \Osheaf_X v_1 \oplus \Osheaf_X v_2 \subset T_X$
is not integrable and $T_X = V \oplus p_{\PP^1}^* T_{\PP^1}$. 
\end{example}

\section{Proof of theorem \ref{theorembogomolovbetter}}
\label{sectionrationallyconnected}

The strategy of the proof of theorem \ref{theorembogomolovbetter} is in a first step 
to construct a fibre space structure on $X$ and to show in a second step that this fibre space structure
comes from a product structure on $X$. The main technical ingredient is a deep theorem by Bogomolov and McQuillan on the
algebraicity of leaves. 

\begin{theorem} \cite[Thm.0.1.a)]{BM01}, \cite[Thm.1.1]{KST05}
\label{theorembogomolov}
Let $X$ be a projective manifold, and let $V \subset T_X$ be an integrable subbundle. Let \holom{f}{C}{X} be
a curve on $X$ such that $f^* V$ is ample. Then all $V$-leaves meeting $f(C)$ are rationally 
connected closed submanifolds of $X$.
\end{theorem}

\begin{lemma}
\label{lemmafoliationrationally}
Let $X$ be a projective rationally connected manifold such that $T_X=V_1 \oplus V_2$ and $V_1 \subset T_X$ is
an integrable subbundle. 
Then there exists a submersion $X \rightarrow Y$ such that $V_1=T_{X/Y}$.
\end{lemma}

{\bf Proof.}
Since $X$ is rationally connected, there exists a (very free) 
rational curve  \holom{f}{\PP^1}{C} on $X$ such that $f^* T_X$ is ample.
Its quotient $f^* V_1$ is then also ample. 
By theorem \ref{theorembogomolov} this 
implies that the $V_1$-leaves passing through a point of $C$ are rationally
connected (closed) subvarieties of $X$. Now apply corollary \ref{corollaryquotientsubmersion}.
$\square$ 

For the convenience of the reader we state the 
classical Ehresmann theorem which allows to deduce properties of the universal covering
from our structure results.

\begin{theorem} \cite[V.,\S 2,Prop.1 and Thm.3]{CLN85}
\label{theoremehresmann}
Let \holom{\phi}{X}{Y} be a submersion of manifolds with an integrable connection, i.e., an
integrable subbundle
$V \subset T_X$ such that $T_X = V \oplus T_{X/Y}$. 
Suppose furthermore one of the following:
\begin{enumerate}
\item $\phi$ is proper.
\item the restriction of $\phi$ to every $V$-leaf is a (not necessarily finite) \'etale map.  
\end{enumerate}
Then \holom{\phi}{X}{Y} is an analytic fibre bundle
with typical fibre $F$.
More precisely, if $\tilde{Y} \rightarrow Y$ is the universal cover, there is 
a representation
$\holom{\rho}{\pi_1(Y)}{\Aut(F)}$ 
such that $X$ is isomorphic to $(\tilde{Y} \times F)/\pi_1(Y)$.
Denote by $\tilde{F} \rightarrow F$ the universal cover of $F$; then the  map 
$\mu: \tilde{Y} \times \tilde{F} \rightarrow \tilde{Y} \times F \rightarrow (\tilde{Y} \times F)/\pi_1(Y) \simeq X$
is the universal cover of $X$ and 
$\mu^* V \simeq p_{\tilde{Y}}^* T_{\tilde{Y}}$ and 
$\mu^* T_{X/Y} \simeq p_{\tilde{F}}^* T_{\tilde{F}}$.
\end{theorem}

\begin{lemma}
\label{lemmaconnectionrationally}
Let $X$ be a projective manifold that admits a submersion \holom{\phi}{X}{Y} on a projective manifold $Y$.
Suppose furthermore that $\phi$ admits a connection, i.e. a vector bundle
$V \subset T_X$ such that $T_X =V \oplus T_{X/Y}$. If $Y$ is rationally connected,  $V$ is integrable
and $X \simeq Y \times F$, where $F$ is a general fibre of $\phi$. 
\end{lemma}

{\bf Proof.}
If $\dim Y = \rk V = 1$ the result is trivial, so we suppose that $\dim Y >1$.
We will construct a covering family of curves on $X$ such that the general 
member is a smooth rational curve
\holom{f}{\PP^1}{C'} that satisfies
\begin{eqnarray*}
f^* T_{X/Y} &\simeq& \Osheaf_{\PP^1}^{\oplus (\dim X - \dim Y)} \qquad \mbox{and}
\\
f^* \bigwedge^2 V &\simeq& \bigoplus_i \Osheaf_{\PP^1} (a_i) 
\qquad \mbox{with} \quad a_i \geq 1 \qquad \forall \ i \geq 1.
\end{eqnarray*}
Granting the construction for the time being, this implies
\[
f^* \Homsheaf(\bigwedge^2 V, T_X/V) \simeq \bigoplus_{i} \Osheaf_{\PP^1} (-a_i)^{\oplus (\dim X - \dim Y)},
\]
so clearly
\[
H^0 (C', \Homsheaf(\bigwedge^2 V, T_X/V)|_{C'})= H^0 (\PP^1, f^* \Homsheaf(\bigwedge^2 V, T_X/V)) = 0,
\]
which shows 
the integrability of $V$ (see subsection \ref{subsectionfoliations}).
Since $\phi$ admits an integrable connection, it follows from the Ehresmann theorem \ref{theoremehresmann} that 
$\phi$ arises as a representation of the fundamental group of $Y$. 
Since $Y$ is simply connected, this 
implies $X \simeq Y \times F$,
where $F$ is a $\phi$-fibre.

We come to the construction of the family.
If $Y$ has dimension at least 3, there
exists a covering family of curves on $Y$ such that the general member is a very free
{\it smooth} rational curve $\holom{f}{\PP^1}{C}$ (cf. \cite[Thm.3.14]{Ko96}). If $Y$ has dimension 2
its relative minimal model is $\PP^2$ or a Hirzebruch surface, so there exists a covering family of free {\it smooth} rational curves.
In both cases $f^* \wedge^2 T_Y \simeq \bigoplus_i \Osheaf_{\PP^1} (a_i)$ with all $a_i$ positive.

Let $\holom{f}{\PP^1}{C}$ be a general member of the covering family. Since $C$ is smooth and $\phi$ is a 
submersion, 
the fibre product $Z:= X \times_Y \PP^1$ is smooth and admits a submersion
\holom{\tilde{\phi}}{Z}{\PP^1}.
Denote by \holom{\mu}{Z}{X} the natural projection; then $\mu^* T_{X/Y}= T_{Z/\PP^1}$
and $T_Z \subset \mu^* T_X$ is a subbundle.
We  consider the sequence of sheaf homomorphisms on $Z$
\[
T_{Z/\PP^1} \hookrightarrow T_Z \hookrightarrow \mu^* T_X \rightarrow \mu^* T_{X/Y}.
\]
The first two maps are the canonical embeddings, while the last one is the projection along $\mu^* V$.
Since $T_{Z/\PP^1} = \mu^* T_{X/Y}$
the map $T_Z \rightarrow \mu^* T_{X/Y}$ has maximal rank at every point.
It follows that $L:=\mu^* V \cap T_Z= \ker(T_Z \rightarrow \mu^* T_{X/Y})$ is a rank 1 subbundle of $T_Z$ such that
\[
T_Z := L \oplus T_{Z/\PP^1}
\]
Since $L$ has rank 1 it is integrable. This shows that \holom{\tilde{\phi}}{Z}{\PP^1} admits an 
integrable connection, so by the Ehresmann theorem $Z \simeq \PP^1 \times F$ where $F$ is a fibre.
It follows that for any $a \in F$, we obtain a rational curve \holom{\mu'}{\PP^1 \times \{a\}}{C'}  
on $X$ such that $\mu'^* T_{X/Y}$ is trivial
and $\mu'^* \wedge^2 V = \oplus_i \Osheaf_{\PP^1} (a_i)$. 
Since the rational curves we use cover a dense open subset on $Y$, the constructed curves cover a dense open subset on $X$. 
$\square$

{\bf Proof of theorem \ref{theorembogomolovbetter}.} 
Lemma \ref{lemmafoliationrationally} yields the existence
of a submersion $X \rightarrow Y$ with connection one of the direct factors $V_j$. Since $X$ is rationally connected,
the manifold $Y$ is rationally connected, so lemma \ref{lemmaconnectionrationally} applies. $\square$

The following lemma will be useful in the next section.

\begin{lemma}
\label{lemmaconnectionbundle}
Let $X$ be a projective manifold that admits a submersion \holom{\phi}{X}{Y} 
with a connection, i.e., a vector bundle
$V \subset T_X$ such that $T_X =V \oplus T_{X/Y}$. Then $\phi$ is an analytic fibre bundle.
\end{lemma}

{\bf Proof.}
In general $V$ is not integrable,
but if $C \subset Y$ is a smooth curve,  the restriction \holom{\phi|_{\fibre{\phi}{C}}}{\fibre{\phi}{C}}{C} is
a smooth map over a curve and $V \cap T_{\fibre{\phi}{C}}$ is a rank 1 bundle that provides a connection
(cf. the proof of lemma \ref{lemmaconnectionrationally} above for details). Since the connection has rank 1 it is integrable,
so $\phi|_{\fibre{\phi}{C}}$ is an analytic bundle. In particular its fibres are isomorphic. Since we can connect
any two points in $Y$ by a chain of smooth curves, this shows that all fibres are isomorphic complex manifolds. 
By the Grauert-Fischer theorem \cite[p.89]{FG65}
this shows that $\phi$ is a fibre bundle. $\square$

\section{Classification of elementary contractions}
\label{sectionuniruled}

This section provides the main technical tools for the proof of theorem \ref{theoremdim4uniruled}. 
Starting with an elementary observation
(lemma \ref{lemmaungenericpositionmorifibrespace}) we classify the structure of elementary Mori contractions of fibre type (proposition 
\ref{propositionfibretype}) in dimension 4. One should note that the arguments employed also hold
in higher dimensions. 
The classification of elementary contractions of birational type (proposition \ref{propositionmoribirational}) 
depends heavily on previous classification
results and therefore does not generalise to higher dimensions.

\begin{lemma}
\label{lemmaungenericpositionmorifibrespace}
Let $X$ be a projective manifold with $T_X = V_1 \oplus V_2$. Let \holom{\phi}{X}{Y} be an elementary
extremal contraction of fibre type. Let $F$ be a general fibre, then, 
after possible renumbering, $T_{F} \subset V_1|_F$. 

If furthermore 
$\dim X - \dim Y = \rk V_1$ or $\dim X - \dim Y +1 = \rk V_1$, the bundle $V_1$ is integrable.
\end{lemma}

{\bf Proof.}
The general fibre $F$ is a Fano manifold,
so there exists
a covering family such that the general member \holom{f}{\PP^1}{F} is a  very free rational curve in $F$, that is
\[
f^* T_F = \oplus_i \Osheaf_{\PP^1} (a_i) 
\qquad \mbox{with} \quad a_i \geq 1 \qquad \forall \ i \geq 1.
\]
We use lemma \ref{lemmasplit}  to obtain
\[
f^* T_X = f^* T_F \oplus \Osheaf_{\PP^1}^ {\oplus \dim Y} = \oplus_i \Osheaf_{\PP^1} (a_i) \oplus \Osheaf_{\PP^1}^ {\oplus \dim Y}.
\]
By the remark after proposition \ref{propositionCP02Prop18} we may assume after possible renumbering that 
$\det V_1$ is $\phi$-ample and $\det V_2$ is $\phi$-trivial.
This implies
$f^* V_2 = \Osheaf_{\PP^1}^{\oplus \rk V_2}$. We denote by \holom{\delta}{T_F}{V_2|_F} the projection on $V_2|_F$ along $V_1|_F$.
Since $f^* T_F$ is ample and $f^* V_2$ is trivial, the restriction of  $\delta$ to $f(\PP^1)$ is zero. 
Since the very free curves cover a dense set in $F$ 
we see that $\delta$ is zero. This is equivalent to $T_F \subset V_1|_F$. 

Suppose now that $\dim X - \dim Y = \rk V_1$ or $\dim X - \dim Y +1 = \rk V_1$. Then we have
$f^* V_1 = \oplus_i \Osheaf_{\PP^1}(a_i)$ or $f^* V_1 = \oplus_i \Osheaf_{\PP^1}(a_i) \oplus \Osheaf_{\PP^1}$.
In both cases $\wedge^2 V_1$ is ample, and $f^* V_2=f^* T_X/V_1$ is trivial, so we obtain 
\[
H^0(\PP^1, \Homsheaf(f^* \wedge^2 V_1, f^* T_X/V_1))=0.
\]
We saw in subsection \ref{subsectionfoliations} that this implies the integrability of $V_1$. 
$\square$

\begin{proposition}
\label{propositionfibretype}
Let $X$ be a smooth projective fourfold with $T_X = V_1 \oplus V_2$ and $\rk V_1=2$. 
Let \holom{\phi}{X}{Y} be an elementary extremal contraction of fibre type. One of the following holds:
\begin{enumerate}
\item $\dim Y=2$. Then \holom{\phi}{X}{Y} is an analytic fibre bundle such that up to renumbering $T_{X/Y}= V_1$.
\item $\dim Y=3$. Then \holom{\phi}{X}{Y} is a $\PP^1$-bundle or conic bundle such that, up to renumbering,
$T_{F} \subset V_1|_F$ for a general fibre $F$. 
\end{enumerate}
\end{proposition}

{\bf Proof.}
By proposition \ref{propositionCP02Prop18} 
we may suppose up to renumbering that $\det V_1$ is $\phi$-ample and $\det V_2$ is $\phi$-trivial.
In particular all fibres have dimension at most $\rk V_1=2$, so $\dim Y \geq 2$. By lemma 
\ref{lemmaungenericpositionmorifibrespace} we have $T_F \subset V_1|_F$ for a general fibre $F$ and the bundle $V_1$ is integrable.

Suppose that $\dim Y=2$. A general fibre $F$ satisfies $T_F =V_1|_F$.
Since the general fibre is a Fano manifold and
a $V_1$-leaf, corollary \ref{corollaryquotientcodim0}  
shows that $\phi$ is smooth. The vector bundle $V_2$ provides a connection, so
by lemma \ref{lemmaconnectionbundle} we know that $\phi$ is a fibre bundle.

Suppose that $\dim Y=3$. For a point $x \in X$, denote by $V_1^x$ the unique $V_1$-leaf containing $x$.
By Ando's result (\cite[Thm.3.1]{An85}) it is sufficient to show that $\phi$ is equidimensional.
We argue by contradiction and suppose that there exists a fibre that has an irreducible component of dimension $2$. 

1st case. There exists a 2-dimensional irreducible component $F$ of a $\phi$-fibre
that is a $V_1$-leaf.

\noindent This implies that the variety $F_{\red}$ is smooth and the corresponding leaf is compact. Since 
$T_{F_{\red}}=V_1|_{F_{\red}}$ and $\det V_1$ is $\phi$-ample, we see that $F$ is a Fano manifold. 
Corollary \ref{corollaryquotientsubmersion} shows that there
exists a submersion \holom{f}{X}{Z} such that all fibres are $V_1$-leaves. 
Since $\phi$ contracts $F$, we can apply the rigidity lemma  \cite[Lemma 1.15]{Deb01} to obtain a dominant factorisation 
\merom{g}{Z}{Y}. 
But this is impossible, since $2 = \dim Z < \dim Y = 3$.

2nd case. For every 2-dimensional irreducible component $F$ of a $\phi$-fibre and every $x \in F$,  the 
set-theoretic intersection $F \cap V_1^x$ is strictly contained in $F$.

\noindent Let $Y^* \subset Y$ be the $\phi$-smooth locus which we 
consider as embedded in the Chow scheme $\chow{X}$. Denote
by $\bar{Y}$ the closure of $Y^*$ in $\chow{X}$, endowed with the reduced structure. 
Let $\Gamma$ be the reduction of the graph over $\bar{Y}$, and \holom{p_{\bar{Y}}}{\Gamma}{\bar{Y}}
and \holom{p_X}{\Gamma}{X} the projections. 
Since every $p_{\bar{Y}}$-fibre is contracted by $\phi$ (this depends only on the homology class), 
the rigidity lemma \cite[Lemma 1.15]{Deb01} implies
the existence of a factorisation \holom{g}{\bar{Y}}{Y}, so that 
we obtain a commutative diagram
\[
\xymatrix{
\Gamma \ar[d]_{p_{\bar{Y}}} \ar[r]^{p_X} & X \ar[d]_\phi
\\
\bar{Y} \ar[r]^g & Y
}
\]
We show that $p_X$ is an isomorphism. Since $p_X$ is birational and $X$ normal
the fibre $\fibre{p_X}{x}$ is connected for every $x \in X$ (Zariski Main Theorem). 
Suppose that it has positive dimension; then there exists a curve $\Delta \subset \bar{Y}$ such 
that for every $y \in \Delta$, we have $x \in p_X(\fibre{p_{\bar{Y}}}{y})$.
Consider now the foliation induced by $p_X^* V_1$ on $\Gamma \subset \bar{Y} \times X$. Since 
a general
$p_{\bar{Y}}$-fibre is contained in a $p_X^* V_1$-leaf and this is a closed condition, every fibre
$\fibre{p_{\bar{Y}}}{y}$ is contained in a $p_X^* V_1$-leaf. Since for $y \in \Delta$ we have
$x \in p_X(\fibre{p_{\bar{Y}}}{y})$ the sets $p_X(\fibre{p_{\bar{Y}}}{y})$ are contained 
in the same $V_1$-leaf $V_1^x$. Furthermore they are contained in the $\phi$-fibre \fibre{\phi}{\phi(x)}.
This shows that the surface $p_X(\fibre{p_{\bar{Y}}}{\Delta})$ is contained in the intersection of $V_1^x$ 
and $\fibre{\phi}{\phi(x)}$.
So there exists a 2-dimensional irreducible component $F$ of \fibre{\phi}{\phi(x)} 
that is equal to $V_1^x$, a contradiction.
This shows that for every $x \in X$ a point, the fibre $\fibre{p_X}{x}$ is a point, so $p_X$ is bijective. 
Since $X$ is smooth and $\Gamma$ reduced this shows that $p_X$ is an isomorphism.

To simplify the notation we identify $\Gamma$ and $X$ via the isomorphism $p_X$.
Let now $F$ be a 2-dimensional component of a fibre. 
The surface $F$ is not contracted by the equidimensional map $p_{\bar{Y}}$.
Let $C \subset F$ be a curve that is not contracted by $p_{\bar{Y}}$.
Since $C \subset F$, it is contracted by $\phi$. Since $\phi$ is
elementary the homology class $[C]$ is a multiple of $[D]$, 
where $D$ is a general $\phi$-fibre. Since the class $(p_{\bar{Y}})_*[D]$ is zero,
the homology class $(p_{\bar{Y}})_* [C]$ is zero. So $p_{\bar{Y}}$ contracts $C$, a contradiction.
$\square$

\vspace{0.5ex}
{\bf Notation.}
Let \holom{\phi}{X}{Y} be a fibration between quasiprojective manifolds. 
The canonical map $\phi^* \Omega_Y \rightarrow \Omega_X$ induces 
a generically surjective sheaf homomorphism 
\holom{T \phi}{T_X}{\phi^* T_Y}. In particular for $\Ssheaf \subset T_X$ a quasicoherent subsheaf,
we obtain a quasicoherent subsheaf $\phi_* (T \phi (\Ssheaf)) \subset T_Y$. 
For a point $y \in Y$ and a point $x \in \fibre{\phi}{y}$ we denote by \holom{(T \phi)_x}{T_{X,x}}{T_{Y,y}}
the tangent map between the vector spaces $T_{X,x}$ and $T_{Y,y}$.

\begin{lemma}
\label{lemmapushdown}
Let $X$ be a quasiprojective manifold with split 
tangent bundle $T_X = V_1 \oplus V_2$.
Let \holom{\phi}{X}{Y} be a fibration
onto a quasiprojective manifold $Y$ such that for a general fibre $F$,
we have 
\[
T_{F} = (V_1 \cap T_F) \oplus (V_2 \cap T_F).
\] 
Suppose furthermore that there exists a Zariski open set $Y^* \subset Y$ such that $Y \setminus Y^*$ 
has codimension at least 2, and such that for 
$y \in Y^*$, the following condition
is satisfied:
\[
 \exists \ x \in \fibre{\phi}{y} \ \mbox{such that} \ \rk (\holom{(T \phi)_x}{T_{X,x}}{T_{Y,y}}) = \dim Y.
\]
For $j=1,2$, the reflexive sheaf $W_j:= (\phi_* (T \phi(V_j)))^{**} \subset T_Y$ 
is a subbundle of $T_Y$ such that $T_Y=W_1 \oplus W_2$.
\end{lemma}

{\bf Proof.}
Step 1. Suppose that $\phi$ is smooth. Then the map \holom{T \phi}{T_X}{\phi^* T_Y} is surjective.
Its restriction to the subbundle $V_j$, denoted by \holom{q_j}{V_j}{\phi^* T_Y} is a morphism of sheaves.
Since $T_{X/Y}= \ker (T \phi)$, we have $\ker q_j= T_{X/Y} \cap V_j$, so the 
rank of $q_j$ at a general point is $\rk V_j - \rk (T_{F} \cap V_j|_F)$. By hypothesis this implies
$\rk q_1 + \rk q_2= \dim Y$. Since $T \phi=q_1 \oplus q_2$ and $T \phi$ has rank equal to $\dim Y$ at every
point, this implies that $q_j$ is a morphism of vector bundles and $\im(q_1) \oplus \im(q_2) = \phi^* T_Y$.
We verify that this induces a splitting of $T_Y$: 
for every fibre $F$, we have 
\[
\im(q_1)|_F \oplus \im(q_2)|_F = \phi^* T_Y|_F  \simeq \Osheaf_F^{\oplus \dim Y},
\] 
so it is elementary to see that $\im(q_j)|_F$ is trivial. This implies $\im(q_j)=\phi^* E_j$ where $E_j$ is a vector
bundle on $Y$, so the splitting pushes down to $Y$.

Step 2. We show the general case.
Let $W_j:= (\phi_*(T\phi(V_j)))^{**} \subset T_Y$; 
then $W_j$ is a reflexive sheaf, so the locus 
$Y' \subset Y$ where $W_1$ and $W_2$ are locally free satisfies $\codim_Y ( Y \setminus Y') \geq 2$. 
Apply the first step to the $\phi$-smooth locus 
to see that $\rk W_1 + \rk W_2 = \dim Y$. 

For $y \in Y' \cap Y^ *$, let $x \in \fibre{\phi}{y}$ be such that the rank of
$\holom{(T \phi)_x}{T_{X,x}}{T_{Y,y}}$ is maximal.
Denote by $V_{j,x} \subset T_{X,x}$ and $W_{j,y} \subset T_{Y,y}$ the subspaces induced by the subbundles $V_j$ and $W_j$. 
Since $(T \phi)_x (V_{j,x})  \subset W_{j,y}$, we have 
\[
T_{Y,y} = \im (T \phi)_x \subset (T \phi)_x(V_{1,x}) +  (T \phi)_x(V_{2,x}) \subset W_{1,y} + W_{2,y} \subset T_{Y,y}.
\] 
Since $\rk W_1 + \rk W_2= \dim Y$ this implies that $W_{1,y} \oplus  W_{2,y} = T_{Y,y}$, 
hence $T_{Y' \cap Y^ *}=W_1|_{Y' \cap Y^ *} 
\oplus W_2|_{Y' \cap Y^ *}$.
Since $Y \setminus (Y^* \cap Y')$ has codimension 2 and $Y$ is smooth, 
we have
\[
T_Y = W_1 \oplus W_2.
\]
In particular the sheaves $W_j$ are locally free.
$\square$

\begin{proposition}
\label{propositionfibretype1}
Let $X$ be a projective fourfold such that $T_X = V_1 \oplus V_2$ with $\rk V_1=\rk V_2 =2$. 
Suppose furthermore that $X$ has the structure of a $\PP^1$-bundle or conic bundle 
\holom{\phi}{X}{Y} such that for a general fibre $F$, we have $T_{F} \subset V_1|_F$.
The map $\phi$ induces a splitting $T_Y = (\phi_*(T \phi (V_1)))^{**} \oplus (\phi_*(T \phi (V_2)))^{**}$.
If $F$ is a singular fibre,  $F$  is isomorphic
to two smooth rational curves intersecting transversally in one point.

Furthermore the $\phi$-singular locus $\Delta \subset Y$ is smooth and satisfies
\[
T_\Delta = (\phi_*(T \phi (V_2)))^{**}|_\Delta.
\] 
\end{proposition}

{\bf Proof.}
If $\phi$ is a $\PP^1$-bundle, the existence of the splitting of $T_Y$ 
follows from lemma \ref{lemmapushdown}.

If $\phi$ is a conic bundle,  a singular fibre is either a double line or two smooth rational curves
intersecting transversally.
The subvariety $D \subset Y$
such that the fibre over every $y \in D$ is a double line has codimension at least 2 
\cite[Prop.1.8,5]{Sar82}. So for $y \in (Y \setminus D)$, there exists a point $x \in \fibre{\phi}{y}$ such that
the tangent map \holom{T \phi}{T_{X,x}}{T_{Y,y}} is surjective.
By lemma \ref{lemmapushdown}, this shows that $\phi$ induces the splitting 
\[
T_Y = (\phi_* (T \phi (V_1)))^{**} \oplus (\phi_* (T \phi (V_2)))^{**}.
\]
In particular $L_1:=(\phi_* (T \phi (V_1)))^{**}$ is a line bundle.

Suppose now that $\phi$ admits a double line $F=\fibre{\phi}{y}$ as a fibre.
The foliation $L_1$ induces around $y$ a germ of a smooth curve $C$ on $Y$ and the foliation $V_1$ induces
around $F_{\red}$ a germ of a smooth surface $S$ on $X$ such that $\phi$ induces a morphism 
\holom{\phi|_S}{S}{C}. The restriction of $N_{F/S}$ to $F_{\red}$ is
a non-trivial torsion line bundle. This contradicts the fact that $F_{\red} \simeq \PP^1$ is simply connected.

The smoothness of the $\phi$-singular locus follows from \cite[Prop.1.8]{Sar82}. Let $C$ be an irreducible
component of a singular fibre, then 
\[
T_X|_C \simeq \sO_{\PP^1} (2) \oplus \sO_{\PP^1}^{\oplus 2} \oplus \sO_{\PP^1} (-1).
\] 
By proposition \ref{propositionCP02Prop18} we know that 
$\det V_1$ is $\phi$-ample and $\det V_2$ is $\phi$-trivial, 
so $T_C \subset V_1|_C$ and $V_2|_C \simeq \sO_{\PP^1}^{\oplus 2}$. 
Let $D \subset X$ be an irreducible effective divisor such that
$\phi(D) \subset \Delta$.  
The canonical map \holom{\gamma}{V_2|_D}{N_{D/X}} is zero
since its restriction to every component $C \subset D$ of a fibre is given by 
\holom{\gamma|_C}{\sO_{\PP^1}^{\oplus 2}}{\sO_{\PP^1}(-1)}, which is the zero map.  
Hence we obtain a splitting of the tangent bundle of the nonsingular locus $D_{\nons} \subset D$ as
\[
T_{D_{\nons}}= (V_1|_{D_{\nons}} \cap T_{D_{\nons}}) \oplus V_2|_{D_{\nons}}.
\]
The inclusion $T_C \subset V_1|_C$ implies
$T_{D_{\nons}/\phi(D_{\nons})}= V_1|_{D_{\nons}} \cap T_{D_{\nons}}$, 
so $\phi(D)=\phi(D_{\nons})$ implies $T_\Delta|_{\phi(D)} = T_{\phi(D)}=W_2|_{\phi(D_{\nons})}$. 
$\square$

\vspace{0.5ex}
{\bf Remark.}
Since elementary contractions \holom{\phi}{S}{C} from a projective surface to a curve are always 
$\PP^1$-bundles, one might expect that in the situation above there are no singular fibres at all.
The following counterexample, for which we thank M. Brunella, shows that this is not true, thereby
correcting and completing \cite[Thm.2.8]{CP02}.
 
Let $S' := \PP^1 \times \PP^1$, we identify $\PP^1 = \C \cup \{ \infty \}$ 
and choose coordinates $(z,w)$ on $S'$. The following map is an
involution on $S'':=S' \setminus \{ (0,0) , (\infty,\infty) \}$.
\[
\holom{\phi'}{S''}{S''} \quad (z,w) \mapsto (z,\frac{z}{w})
\]
If $C_z:=\{z\} \times \PP^1$ is a fibre of the projection $pr_1$ on the first factor,
then $\phi'|_{C_z}$ is the involution on $\PP^1$ with fixed points $\sqrt{z}$ and
$- \sqrt{z}$.
Blow up $S'$ in $(0,0)$ and $(\infty,\infty)$ to resolve the indeterminacies of
$\phi'$. If we denote by $\holom{\mu}{S}{S'}$ the blow-up, $\phi'$ lifts to
an involution $\phi$ of $S$. The total transform of $\{0\} \times \PP^1$ has two
irreducible components $E_1$ and $E_2$ such that $\phi(E_1)=E_2$ and $\phi(E_2)=E_1$
(an analogous statement holds for the total transform of $\{\infty\} \times \PP^1$).

Let $T:=\C/(\Z \oplus i \Z)$ be an elliptic curve and $\holom{\psi}{T}{T} \quad t \mapsto t+\frac{1}{2}$. Then
$T' \simeq T/\{\id_T,\psi\}$ is an elliptic curve. 
Define $X:= (S \times T)/ \{\id_S \times \id_T, \phi \times \psi\}$, then $X$ is a smooth
projective variety with split tangent bundle. The map $(pr_1 \circ \mu) \times \id_T$ induces
a morphism 
\[
\holom{f}{X}{\PP^1 \times T'}.
\]
Let us show that $f$ is an elementary contraction. 
Let  $F_0' = E_1' + E_2'$ be a singular fibre and let $F_0 = E_1 + E_2$ be its lifting to $S \times T$.
Then $F_0 \subset S \times \{t\}$ for some $t \in T$ and
clearly 
\[
E_1 = (\phi \times \psi) (\tilde{E})
\] 
where $\tilde{E}$ is the copy of $E_2$ in $S \times \{t+\frac{1}{2}\}$.
Since $\tilde{E}$ is a deformation of $E_2$ in $S \times T$, this shows that  
$E_2'$ is a deformation
of $E_1'$ in $X$, hence $E_1'=E_2'$ in $N_1(X)$. 
The special fibre $F_0$ is homologous to a general fibre $F$, so
\[
F = F_0 = E_1' + E_2' = 2 E_1'
\] 
in $N_1(X)$. This shows that $f$ is the elementary contraction of the ray generated by $E_1'$.

\vspace{0.5ex}
{\bf Notation.} Let \holom{\phi}{X}{Y} be an elementary contraction of birational type.
If the  exceptional locus $E$ of $\phi$ is irreducible, 
let $k:=\dim E$ and $l:= \dim \phi(E)$. 
The birational contraction is then said to be of type $(k,l)$. 

\begin{proposition}
\label{propositionmoribirational}
Let $X$ be projective fourfold such that $T_X = V_1 \oplus V_2$ with $\rk V_1=\rk V_2=2$,
and let \holom{\phi}{X}{Y} be an elementary contraction of birational type. 
Then $Y$ is smooth and $\phi$ is the blow-up of a smooth $2$-dimensional subvariety of the manifold $Y$.
Set $W_j:= (\phi_* (T\phi(V_j)))^{**}$; then the tangent bundle of $Y$ splits as 
\[
T_Y = W_1 \oplus W_2.
\]
If the universal covering \holom{\mu}{\tilde{Y}}{Y} splits as $\tilde{Y} \simeq Y_1 \times Y_2$  
such that $\mu^* W_j = p_{Y_j}^* T_{Y_j}$,  the analogous statement holds for $X$.
\end{proposition}

{\bf Proof.}
Let $E$ be the exceptional locus of the contraction
and
$F$ an irreducible component of a non-trivial fibre. Then by proposition
\ref{propositionCP02Prop18} we have $\dim F \leq 2$. In particular the contraction
can't be of type $(3,0)$.
Contractions of type $(3,1)$ have been classified by Takagi \cite[p.316]{Ta99}: the map 
\holom{\phi|_E}{E}{\phi(E)} is a $\PP^2$-bundle or quadric bundle, so 
the general fibre is reduced and isomorphic to $\PP^2$, $\PP^1 \times \PP^1$, or a quadric cone.
If $\phi$ is of type $(3,2)$,  a fibre of dimension 2 is isolated, so $F$ is reduced and 
either $\PP^2$,  $\PP^2 \cup \PP^2$, $\PP^1 \times \PP^1$, or a quadric cone \cite[Thm.4.6]{AW96}.
Elementary fourfold contractions such that the exceptional locus is not an irreducible divisor 
have been classified by Kawamata \cite[Thm.1.1]{Ka89}:
the positive dimensional fibres are reduced and isomorphic to $\PP^2$.

We will exclude case by case the existence of 2-dimensional fibres. 
The strategy is to
choose appropriately rational curves \holom{f}{\PP^1}{F} such that $f^* T_F$ is nef
so that we can use equation (\ref{equationfibresplit}) to compute $f^* T_X$.

Case 1. $F \simeq \PP^2$ or $F \simeq \PP^2 \cup \PP^2$.

If $F \simeq \PP^2$, let \holom{f}{\PP^1}{F} be a line $C=f(\PP^1)$, then $f^* N_{F/X}^*$ is nef, so
by equation (\ref{equationfibresplit})
\[
f^* T_X = \Osheaf_{\PP^1} (2) \oplus \Osheaf_{\PP^1} (1) \oplus \Osheaf_{\PP^1} (-a) \oplus \Osheaf_{\PP^1} (-b),
\]
where $a \geq 0$, $b \geq 0$ and $a+b>0$ (otherwise the deformations of $C$ would cover $X$).
Since the curve is of splitting type (\ref{equationsplittingtype}), we have up to renumbering  $\Osheaf_{\PP^1} (2) \subset f^* V_1$.
 Since $C$ is a line in projective space it deforms keeping a point fixed.
This shows that $\Osheaf_{\PP^1}(2) \subset f^* V_1$ implies $T_{\PP^2}=V_1|_{\PP^2}$. In particular 
$f^* V_1 = \Osheaf_{\PP^1} (2) \oplus \Osheaf_{\PP^1} (1)$, so $f^* V_2 = \Osheaf_{\PP^1} (-a) \oplus \Osheaf_{\PP^1}(-b)$. 
This implies that
$f^* \det V_2 = \Osheaf_{\PP^1}(-a-b)$ is not $\phi$-trivial, a contradiction to proposition \ref{propositionCP02Prop18}.
The same argument works in the case
$F \simeq \PP^2 \cup \PP^2$.

Case 2. $F \simeq \PP^1 \times \PP^1$.

Choose a ruling line \holom{f}{\PP^1}{F} with image $C=f(\PP^1)$. Then
\[
f^* T_X = \Osheaf_{\PP^1} (2) \oplus \Osheaf_{\PP^1} \oplus 
\Osheaf_{\PP^1} (-a) \oplus \Osheaf_{\PP^1} (-b)
\]
where $a \geq 0$ and $b \geq 0$ and $a+b>0$. Up to renumbering, this implies $\Osheaf_{\PP^1}(2) \subset f^* V_1$.
By  proposition \ref{propositionCP02Prop18} we see that
$\det V_1$ is $\phi$-ample and $\det V_2$ is $\phi$-trivial.
Now choose a line \holom{f'}{\PP^1}{F} from
the second ruling that
is transversal to the first one. Then
\[
f'^* T_X = \Osheaf_{\PP^1} (2) \oplus \Osheaf_{\PP^1} \oplus 
\Osheaf_{\PP^1} (-a) \oplus \Osheaf_{\PP^1} (-b)
\]
with coefficients as above.
Then again $\Osheaf_{\PP^1}(2) \subset f'^* V_1$, since otherwise $f'^* (c_1(V_1)) = 0$ (for details 
cf. the proof of \cite[Lemma 1.3]{CP02}), which would contradict
the $\phi$-ampleness of $\det V_1$. Since the lines are transversal, we obtain
\[
\Osheaf_{\PP^1}(2) \oplus \Osheaf_{\PP^1} = f^* T_F = f^* V_1.
\]
As in the first case we obtain $f^* V_2 = \Osheaf_{\PP^1} (-a) \oplus \Osheaf_{\PP^1}(-b)$. Clearly
$f^* \det V_2 = \Osheaf_{\PP^1}(-a-b)$ is not trivial, a contradiction. 

Case 3. $F$ is a quadric cone.

Let \holom{f}{\PP^1}{F} be a line passing through the vertex $x$ of the cone $F$. Campana and Peternell
have shown in \cite[Thm.3.6]{CP02} that $C=f(\PP^1)$ is of splitting type (\ref{equationsplittingtype}), 
so up to renumbering, we have $\Osheaf_{\PP^1}(2) \subset f^* V_1$.
In particular $T_{C,x} \subset V_{1,x}$ for every such line. Since $x$ is a singularity of $F$ and the vector spaces
$T_{C,x}$ generate the Zariski tangent space of $F$ in $x$, they generate a subspace of dimension at least 3 in 
$V_{1,x}$. This contradicts $\rk V_1=2$.

We summarize: the morphism $\phi$ is a birational contraction of divisorial type such that all non-trivial
fibres are of pure dimension 1. So the structure of $\phi$ follows from a result
by Ando \cite[Thm.2.1]{An85}. 

Let $F$ be a positive dimensional fibre;
then $T_X|_F \simeq \sO_{\PP^1} (2) \oplus \sO_{\PP^1}^{\oplus 2} \oplus \sO_{\PP^1} (-1)$. 
Up to renumbering $\det V_1$ is $\phi$-ample and $\det V_2$ is $\phi$-trivial, 
so $T_F \subset V_1|_F$ and $V_2|_F \simeq \sO_{\PP^1}^{\oplus 2}$. 
The splitting of the tangent bundle $T_Y$ is  a consequence of lemma \ref{lemmapushdown}.
Let $E$ be the exceptional divisor; then the canonical map \holom{\gamma}{V_2|_E}{N_{E/X}} is zero
since its restriction to the positive dimensional fibres is given by
\holom{\gamma|_F}{\sO_{\PP^1}^{\oplus 2}}{\sO_{\PP^1}(-1)}, which is the zero map.  
It follows easily that $T_E= (V_1|_E \cap T_E) \oplus V_2|_E$ and $T_F \subset V_1|_F$ implies
$T_{E/\phi(E)}= V_1|_E \cap T_E$, so we have $T_{\phi(E)}=W_2|_{\phi(E)}$.

Let now \holom{\mu}{\tilde{Y}}{Y} be the universal covering map and
suppose that $\tilde{Y} \simeq Y_1 \times Y_2$  
and $\mu^* W_j = p_{Y_j}^* T_{Y_j}$.
Note that this implies the integrability of $W_1$ and $W_2$.
Since $\pi_1(X)=\pi_1(Y)$, the pull-back $\tilde{X} = X \times_Y \tilde{Y}$
is the universal covering of $X$.
We have seen that $\phi$ is the blow-up of $Y$ along the $W_2$-leaf $\phi(E)$.
It follows that $\tilde{\phi}$ is the blow-up of $\tilde{Y}$ along a leaf of the foliation 
$\mu^* W_2=p_{Y_2}^* T_{Y_2}$.
So there exists a $y \in Y$ such that 
\[
\tilde{X} \simeq \Bl{y \times Y_2}{\tilde{Y}} \simeq \Bl{y}{Y_1} \times Y_2 =: X_1 \times X_2,
\]
where $\Bl{A}{B}$ denotes the blow-up of a manifold $B$ along a submanifold $A$. 
Since 
\[
\tilde{\mu}^* V_j = \tilde{\phi}^* \circ \mu^* W_j = \tilde{\phi}^* p_{Y_j}^* T_{Y_j},
\]
we have $\tilde{\mu}^* V_1=p_{X_1}^* T_{X_1}$ and $\tilde{\mu}^* V_2=p_{X_2}^* T_{X_2}$. 
$\square$

\section{Proof of theorem \ref{theoremdim4uniruled}}

{\bf Proof.}
By proposition \ref{propositionmoribirational} the splitting of $T_X$ is stable
under birational contractions. By the general minimal model program it is clear that after finitely many 
birational contractions
$X = X^0 \rightarrow X^1 \rightarrow \ldots \rightarrow X^n$ 
we obtain either a variety $X^n$ with nef canonical bundle or an elementary contraction of fibre type.
Since $X$ is uniruled, the variety $X^n$ is uniruled so we are in the second case.

By proposition \ref{propositionmoribirational} it is sufficient to
show conjecture \ref{conjecturebeauville} for $X^n$. To simplify the notation we suppose that $X=X^n$.
Suppose now that $V_1$ and $V_2$ are integrable. 
By  proposition \ref{propositionfibretype} there are the following two cases.

1st case. The variety $X$ is an analytic fibre bundle  
\holom{\phi}{X}{Y} such that up to renumbering $T_{X/Y}=V_1$.
Since $V_2$ is integrable, we conclude with the Ehresmann theorem.

2nd case. The variety $X$ is a $\PP^1$-bundle or conic bundle  
\holom{\phi}{X}{Y} such that up to renumbering $T_{F} \subset V_1|_F$ for a general fibre.
By proposition \ref{propositionfibretype1} we have $T_Y = (\phi_* (T \phi (V_1)))^{**} \oplus (\phi_* (T \phi (V_2)))^{**}$
where $L:=(\phi_* (T \phi (V_1)))^{**}$ is a line bundle.
By theorem \ref{theorembrunella} applied to $Y$, there are two
subcases.

Case a. The variety $Y$ is a $\PP^1$-bundle \holom{\psi}{Y}{Z} such that $T_{Y/Z}=L$. 
Then \holom{\psi \circ \phi}{X}{Z} is a proper submersion
with integrable connection $V_2$ and we conclude with the Ehresmann theorem \ref{theoremehresmann}.

Case b. The bundle $V:=(\phi_* (T \phi (V_2)))^{**}$ is integrable, 
so the universal covering \holom{\mu}{\tilde{Y}}{Y} satisfies $\tilde{Y} \simeq S \times C$
such that $\mu^* V = p_S^* T_S$ and $\mu^* L = p_C^* T_C$. Furthermore we have a commutative diagram
\[
\xymatrix{
X' \ar@/_/[dd]_q \ar[r]^{\tilde{\mu}} \ar[d]^{\tilde{\phi}} & X  \ar[d]^\phi
\\
\tilde{Y}  \ar[d]^{p_S} \ar[r]^\mu & Y
\\
S 
}
\]
where \holom{\tilde{\mu}}{X':= X \times_Y \tilde{Y}}{X} is \'etale. 
By proposition \ref{propositionfibretype1} the morphism $\tilde{\phi}$ does not have any multiple fibres, so the 
almost smooth morphism 
\holom{q:=p_S \circ \tilde{\phi}}{X'}{S} is a submersion 
such that $\tilde{\mu}^* V_1 = T_{X'/S}$. 
The morphism $q$ is not necessarily proper, so to apply
theorem \ref{theoremehresmann} we have to verify that the restriction of $q$ to
every $\tilde{\mu}^*  V_2$-leaf is a covering map.
Let $\mathfrak{V}$ be a $\tilde{\mu}^*  V_2$-leaf, 
then 
\[
(\tilde{\phi}_* (T \tilde{\phi} (\tilde{\mu}^* V_2)))^{**} = \mu^* V = p_S^* T_S
\]
implies that $\phi(\mathfrak{V}) = S \times c$ for some $c \in C$. 
Set $Z:=\fibre{\tilde{\phi}}{S \times c}$, then \holom{\tilde{\phi}|_Z}{Z}{S \times c}
is proper. By proposition \ref{propositionfibretype1} we know that $S \times c$ is either disjoint from
the $\tilde{\phi}$-singular locus or contained in it.
In the first case 
apply \cite[V.,\S 2,Prop.1]{CLN85}
to see 
that \holom{\tilde{\phi}|_{\mathfrak{V}}}{\mathfrak{V}}{S \times c} is a covering map.
In the second case 
let \holom{\nu}{Z'}{Z} be the normalisation of $Z$, then 
\holom{\tilde{\phi} \circ \nu}{Z'}{S \times c}
is a fibre bundle with typical fibre two connected components isomorphic to $\PP^1$
and $\nu^* \tilde{\mu}^* V_2 \simeq \nu^* \tilde{\phi}^* T_{S \times c}$. 
Apply again \cite[V.,\S 2,Prop.1]{CLN85} to see that the restriction of 
$\tilde{\phi} \circ \nu$ to every $\nu^* \tilde{\mu}^* V_2$-leaf is a covering map.
This implies
that \holom{\tilde{\phi}|_{\mathfrak{V}}}{\mathfrak{V}}{S \times c} is a covering map
 
Since \holom{p_S|_{S \times c}}{S \times c}{S} is an isomorphism we obtain in both cases
that the restriction of $q$ to $\mathfrak{V}$ is a covering map. 
 $\square$

\end{document}